\documentclass[12pt,twoside]{article}
\usepackage{graphicx,graphics, amsfonts, amsbsy,amssymb, amsmath, cite, array,arydshln}

\allowdisplaybreaks

\newtheorem{theorem}{Theorem}[section]

\newtheorem{problem}{Problem}
\newtheorem{lemma}[theorem]{Lemma}
\newtheorem{corollary}[theorem]{Corollary}

\usepackage{tikz,pgfplots}
\usetikzlibrary{chains,positioning}
\usetikzlibrary{decorations.pathmorphing}

\def\qed{\nolinebreak\hfill\rule{.2cm}{.2cm}\par\addvspace{.5cm}}

\begin{document}
\title{On $A_{\alpha}$-spectrum of joined union of graphs and its applications to power graphs of finite groups}
\author{ Bilal A. Rather$ ^{a} $, Hilal A. Ganie$^{b}$,  S. Pirzada$ ^{c} $  \\
$^{a,c}${\em  Department of Mathematics, University of Kashmir, Srinagar, India}\\
$^{b}$ {\em Department of School Education, JK Govt. Kashmir, India}\\
$^a$bilalahmadrr@gmail.com,~~ $^b$hilahmad1119kt@gmail.com \\
$^c$pirzadasd@kashmiruniversity.ac.in} 
\date{}

\pagestyle{myheadings} \markboth{Bilal, Hilal and Pirzada}{On $ A_{\alpha} $-spectrum of joined union and its applications to power graphs of finite groups}
\maketitle
\vskip 3mm

\noindent{\footnotesize \bf Abstract.}
For a simple graph $G$, the generalized adjacency matrix  $A_{\alpha}(G)$ is defined as   $A_{\alpha}(G)=\alpha D(G)+(1-\alpha)A(G), \alpha\in [0,1]$, where $A(G)$ is the adjacency matrix and $D(G)$ is the diagonal matrix of vertex degrees of $G$. This matrix generalises the spectral theories of the adjacency matrix and the  signless Laplacian matrix of $G$. In this paper, we find $ A_{\alpha} $-spectrum of the joined union of graphs in terms of spectrum of adjacency matrices of its components and the zeros of the characteristic polynomials of an auxiliary matrix determined by the joined union. We determine the $ A_{\alpha}$-spectrum of join of two regular graphs, the join of a regular graph with the union of two regular graphs of distinct degrees. As an applications, we investigate the $ A_{\alpha} $-spectrum of certain power graphs of finite groups.

\vskip 3mm



\section{Introduction}

Let $ G $ be a simple graph with $ n $ vertices and $ m $ edges having vertex $V(G) = \{v_{1},v_{2},\dots,v_{n}\} $ and edge set $ E(G)= \{e_{1},e_{2},\dots,e_{m}\} $. The set of vertices adjacent to  a vertex $v\in V(G)$, denoted by $N(v)$, refers to the \textit{neighborhood} of $v.$ The \textit{degree} of a vertex  $v,$ denoted by $d_{G}(v)$ (we simply write $d_v$) is the cardinality of $N(v)$. A graph is called \textit{regular} if each of its vertices have the same degree. We write $ v_{i}\sim v_{j} $, if $ v_{i} $ and $ v_{j} $ are connected by an edge, otherwise we write $ v_{i}\nsim v_{j}. $  The adjacency matrix $A=(a_{ij})$ of $G$ is a $(0, 1)$-square matrix of order $n$ whose $(i,j)$-entry is equal to $1$, if $v_i\sim v_{j} $ and equal to $0$, otherwise. Let $ D(G)=diag(d_{1},d_{2},\dots,d_{n}) $ be the diagonal matrix of vertex degrees $ d_{i}=d_{G}(v_{i}) , i=1,2,\dots,n$ of  $ G. $ The matrices $ L(G)=A(G)-D(G) $ and $ Q(G)=A(G)+D(G) $ are called the Laplacian matrix and the signless Laplacian matrix, respectively of the graph $G$. Both $ L(G) $ and $ Q(G) $ are positive semidefinite matrices having real eigenvalues so that their eigenvalues can be put as $ \mu_{1}(G)\geq\mu_{2}(G)\geq\dots\geq\mu_{n}(G)=0 $ and $ q_{1}(G)\geq q_{2}(G)\geq\dots\geq q_{n}(G)$, respectively.

Nikiforov \cite{vn} proposed to study the convex combinations $ A_\alpha(G) $ of $ A(G) $ and $ D(G) $ defined by
$ A_{\alpha}(G)=\alpha D(G) +(1-\alpha)A(G)$, $ 0\leq \alpha \leq 1$. Clearly, $ A(G)=A_{0}(G)$, $D(G)=A_{1}(G) $ and $ 2A_{\frac{1}{2}}(G)=D(G)+A(G)=Q(G)$. We further note that
$ A_{\alpha}(G)-A_{\gamma}(G)=(\alpha-\gamma)(D(G)-A(G))=(\alpha-\gamma)L(G)$, where $L(G)$ is the Laplacian matrix of $G$. Therefore, the matrix $A_{\alpha}(G)$ is a generalisation of the adjacency matrix $A(G)$ and the signless Laplacian matrix $Q(G)$ of $G$.  Let $A(G)$ be the adjacency matrix of the graph $G$ and let $B$ be a real diagonal matrix of order $n$. Bapat et al. \cite{bkp} defined the matrix $L^{'}=B-A(G)$ and called it the perturbed Laplacian matrix of $G$. For $\alpha\ne 1$, it is easy to see that 
\begin{align*}
A_{\alpha}(G)=\alpha D(G)+(1-\alpha)A(G)=(\alpha-1)\Big(\frac{\alpha}{\alpha-1}D(G)-A(G)\Big).
\end{align*}  
Clearly $\frac{\alpha}{\alpha-1}D(G)$ is a diagonal matrix with real entries, so that $A_{\alpha}(G)$ is a scaler multiple of a perturbed Laplacian matrix. This is another motivation to study the spectral properties of $A_{\alpha}(G)$. For some recent papers on the spectral properties of $A_{\alpha}(G)$, we refer to  \cite{vn,vng,dl,cw,sbhr} and the references therein.

The matrix $A_{\alpha}(G)$ is a real symmetric matrix, therefore we can arrange its eigenvalues as 
$\lambda_{1}(A_{\alpha}(G))\geq \lambda_{2}(A_{\alpha}(G))\geq \cdots \geq \lambda_{n}(A_{\alpha}(G)),$ where $ \lambda_{1}(A_{\alpha}(G)) $ is called the \textit{generalized adjacency spectral radius} of $G$. For a connected graph $G$, the matrix $A_{\alpha}(G) $ is non-negative and irreducible (for $\alpha\ne 1$).  Therefore, by Perron-Frobenius theorem, $ \lambda(A_{\alpha}(G))$ is the simple (with multiplicity $1$) eigenvalue and there is a unique positive unit eigenvector $X$ corresponding to $ \lambda(A_{\alpha}(G)),$ which is called the \textit{generalized adjacency Perron vector} of $ G.$

As usual $ K_{n},~ K_{1,n-1}~ P_n \text{and}~ C_{n} $ respectively, denote the complete graph, the star, the path and the cycle graph of order $n$. For other notations and terminology we refer to \cite{cds,hj,sp, roman}.

The rest of the paper is organized as follows. In Section 2, we obtain the $ A_{\alpha} $-spectrum of the joined union of regular graphs $G_1,G_2,\dots,G_n$ in terms of their adjacency spectrum and the spectrum of an auxiliary matrix. As a consequence, we obtain $ A_{\alpha} $-spectrum of some well known graphs. In  Section 3, we investigate the $ A_{\alpha} $-spectrum of the power graphs of some finite groups and raise some problems for further work.

\section{$ A_{\alpha} $-spectrum of the joined union of graphs}

Consider an $ n \times n $ matrix
\begin{equation*}
M= \begin{pmatrix}
A_{1,1} & A_{1,2} & \cdots & A_{1,s} \\
A_{2,1} & A_{2,2} & \cdots & A_{2,s} \\
\vdots & \vdots & \cdots & \vdots \\
A_{s,1} & A_{s,2} & \cdots & A_{s,s} \\
\end{pmatrix}, \end{equation*}
whose rows and columns are partitioned according to a partition $
P=\{ P_{1}, P_{2},\dots , P_{m}\} $ of the set $X= \{1,2,\dots,n\}. $ The quotient matrix $ \mathcal{Q} $ (see \cite{BH, cds}) is an $s \times s$ matrix whose entries are the average row sums of the blocks $ A_{i,j} $ of $ M$. The partition $ P $ is said to be {\it equitable} if each block $ A_{i,j} $ of $ M $ has constant row (and column) sum and in this case the matrix $ \mathcal{Q} $ is called as {\it equitable quotient matrix}. A vertex partition $\{ V_{1},V_{2},\dots,V_{s} \}$ of the vertex set $ V(G) $ of $ G $ is equitable if for each $ i $ and for all $ u,v\in V_{i}$, we have $|N(u)\cap V_{j}|=|N(v)\cap V_{j}|$, for all $ j. $ In general, the eigenvalues of $ \mathcal{Q}  $ interlace the eigenvalues of $ M $. In case the partition is equitable, the following lemma  was obtained in \cite{BH}.

\begin{lemma}\cite{BH, cds}
If the partition $ P $ of $ X $ of matrix $ M $ is equitable, then each eigenvalue of $ \mathcal{Q} $ is an eigenvalue of $ M. $
\end{lemma}

\indent Consider a graph $G$ having symmetry so that its associated matrix is written in the form
\begin{align}\label{pe1}
M=\begin{pmatrix}
X & \beta & \cdots &\beta & \beta \\
\beta^T & B & \cdots & C & C\\
\vdots & \vdots & \cdots & \vdots &\vdots\\
\beta^T & C & \cdots &B & C\\
\beta^T & C & \cdots & C & B
\end{pmatrix},\end{align}
where $X\in  R^{t\times t},~\beta \in R^{t\times s}$ and $B,C\in R^{s\times s}$, such that $n=t+ cs$, with $c$ being the number of copies of $B$.   The spectrum of this matrix can be obtained as the union of the spectrum of smaller matrices using the following technique given in \cite{ft}. In the statement of the theorem, $\sigma^{(k)}(Y)$ indicates the multi-set formed by $k$ copies of the spectrum of $Y$, denoted by $\sigma(Y)$.

\begin{lemma}\label{Lemma p2}
Let $M$ be a matrix of the form given in \eqref{pe1}, with $c\geq 1$ copies of the
block $B$. Then
\begin{enumerate}
\item[(i)] $\sigma(B-C)\subseteq \sigma(M)$ with multiplicity $c-1$;
\item[(ii)]  $\sigma(M)\setminus \sigma^{(c-1)}(B-C)=\sigma(M^{\prime})$ is the set of the remaining $t+s$ eigenvalues of $M$, where
$M^{\prime}=\begin{pmatrix}
X & \sqrt{c}\beta  \\
\sqrt{c}\beta^T & B+(c-1)C
\end{pmatrix}.$
\end{enumerate}
\end{lemma}

Let $G(V,E)$ be a graph of order $n$ and $ G_{i}(V_{i}, E_{i})$ be graphs of order $n_i,$ where $i=1,\ldots, n $. The  \textit{joined union} $ G[G_{1},\ldots, G_{n}] $ is the graph $ H(W, F) $ with
\begin{eqnarray*}
W=\bigcup_{i=1}^{n}V_{i}~~~\text{and}~~F=\bigcup_{i=1}^{n}E_{i}\cup\bigcup_{\{v_{i}, v_{j}\}\in E}V_{i}\times V_{j}.
\end{eqnarray*}
In other words, the joined union is the union of graphs $ G_{1},\ldots, G_{n} $ together with the edges $ v_{ik}v_{jl}, ~ v_{ik}\in G_{i} $ and $ v_{jl}\in G_{j} $, whenever $ v_{i}v_{j} $ is an edge in $ G $. Clearly, the usual join of two graphs $ G_1 $ and $G_2$ is a special case of the joined union $ K_{2}[G_1, G_2]=G_{1}\triangledown G_{2} $ where $ K_{2} $ is the complete graph of order $2$. From the definition, it is clear that the diameter of the graph $G[G_{1},\ldots, G_{n}] $ is same as the diameter of the graph $G$, for every choice of the connected graphs  $G_1,G_2,\dots,G_n$, $n\geq 3$ and $G \ncong K_{n}$.

Now, we obtain the $ A_{\alpha} $-spectrum of the joined union of regular graphs $G_1,G_2,\dots,G_n$ in terms of the adjacency spectrum of the graphs $G_1,G_2,\dots,G_n$ and the eigenvalues of an auxiliary matrix.

\begin{theorem} \label{joined union}
Let $G$ be a graph of order $n\ge 3$ having $m$ edges. Let $ G_{i}$ be $r_{i}$-regular graphs of order $ n_{i} $ having adjacency eigenvalues $\lambda_{i1}=r_{i}\geq \lambda_{i2}\geq\ldots\geq \lambda_{in_{i}}, $ where $ i=1,2, \ldots, n$.  The $ A_{\alpha} $-spectrum of the joined union graph $ G[G_{1},\ldots, G_{n}] $ of order $ N=\sum\limits_{i=1}^{n}n_{i} $ consists of the eigenvalues $ \alpha (r_{i}+\alpha_i )+(1-\alpha)\lambda_{ik}(G_{i})$ for $ i=1,\ldots,n $ and $ k=2,3,\ldots, n_{i} $, where $\alpha_i=\sum\limits_{v_j\in N_{G}(v_{i})}n_{i}$ is the sum of the cardinalities of the graphs $G_j, j\ne i$ which corresponds to the neighbours of vertex $ v_{i}\in G $. The remaining $n$ eigenvalues are given by the matrix
\begin{equation*}
M =\begin{pmatrix}
\psi_{11}& (1-\alpha)\psi_{12}&\ldots & (1-\alpha)\psi_{1n}\\
(1-\alpha)\psi_{21}&\psi_{22}& \ldots & (1-\alpha)\psi_{2n}\\
\vdots &\vdots  &\ddots &\vdots\\
(1-\alpha)\psi_{n1}& (1-\alpha)\psi_{n2}&\ldots &\psi_{nn}
\end{pmatrix},
\end{equation*} 
where $\psi_{ii}=\alpha \alpha_i+r_i$, for $i=1,2,\dots,n$ and for $i\ne j$,  $\psi_{ij}=n_j$, if $v_i\sim v_j$, while as  $\psi_{ij}=0$, if $v_i\nsim v_j$. 
\end{theorem}
{\bf Proof.}
Let  $V(G)=\{ v_{1}, \ldots, v_{n}\}$ be the vertex set of $ G $ and let $V(G_i)=\{ v_{i1}, \ldots, v_{in_i}\}$ be the vertex set of $G_i$, for $i=1,2,\dots,n$. Let  $H=G[G_{1}, \ldots, G_{n}]$ be the joined union  of the graphs $G_1,G_2,\dots, G_n$. Clearly, order of $H$ is $N=\sum\limits_{i=1}^{n}n_i$. We first compute the degree of each vertex in $V(H)$. Evidently, the degree of each vertex $v_{ij}\in V(H)$, for $1\le i\leq n$ and $1\leq j\leq n_i$, is the degree inside $ G_{i} $ plus the sum of cardinalities of the copies of $ G_{j}, j\ne i $ which corresponds to the neighbours of the vertex $ v_{i} $ in  $ G $. That is, for  each $v_{ij}\in V(G_{i})$, we have
\begin{equation}\label{eq23}
d_{H}(v_{ij})=r_{i}+\sum\limits_{v_j\in N_{G}(v_i)}n_{j}=r_i+\alpha_i,
\end{equation} where $\alpha_i=\sum\limits_{v_j\in N_{G}(v_i)}n_{j}$. Obviously, $d_{H}(v_{ij})$ is same for each vertex in $G_i$, $1\le j\le n_i$.
Label the vertices in $ H $ from the vertices in $ G_{1} $ to the vertices in $ G_{n} $. With this labelling, the $ A_{\alpha} $-matrix of $H$ can be written as
\begin{equation*}
A_{\alpha}(H)=\begin{pmatrix}
h_1 & (1-\alpha) a(v_1,v_2)& \ldots & (1-\alpha) a(v_1,v_n) \\
(1-\alpha)a(v_2,v_1) & h_2 & \ldots & (1-\alpha) a(v_2,v_n)\\
\vdots &\vdots &\ddots &\vdots\\
(1-\alpha) a(v_n,v_1) & (1-\alpha)a(v_n,v_2) &\ldots & h_n
\end{pmatrix},
\end{equation*}
where,
\begin{equation*}
h_{i}=\alpha (r_{i}+\alpha_i)I_{n_{i}}+(1-\alpha)A(G_{i}), \qquad  \text{for}~ i=1,2,\ldots,n,
\end{equation*}
and $ a(v_{i}, v_{j}) =J_{n_{i}\times n_{j}} $, if $ v_{i}\sim v_{j} $ in $G$ and $ 0_{n_{i}\times n_{j}} $, otherwise.
For $ i=1,2,\ldots,n, A(G_i)$ is the adjacency matrix of $ G_{i}$, $J_{n_i\times n_{j}}$ is the matrix having all entries $ 1 $, $ 0_{n_{i}\times n_{j}} $ is the zero matrix of order $ n_{i}\times n_{j} $ and $ I_{n_i}$ is the identity matrix of order $n_i$. \\
\indent As $ G_{i} $ is an $ r_{i} $ regular graph, so the all one vector $e_{n_i}=(1,1,\dots,1)^T$ with $ n_{i} $ entries is the eigenvector of the adjacency matrix $ A(G_{i}) $ corresponding to the eigenvalue $r_{i} $ and all other eigenvectors are orthogonal to $e_{n_i}.$ Let $ \lambda_{ik}$, $2\leq k\leq n_i$, be any eigenvalue of $ A(G_{i})$ with the corresponding eigenvector $X=(x_{i1},x_{i2},\dots,x_{in_i})^T$ satisfying $e_{n_i}^TX=0.$ Clearly, the column vector $X$ can be regarded as a function defined on $ V(G_i) $ assigning the vertex $ v_{ij} $ to $ x_{ij} $, that is, $ X(v_{ij})=x_{ij} $ for $ i=1,2,\ldots,n $ and $j=1,2,\dots,n_i$. Now, consider the vector $Y=(y_{1},y_{2},\dots,y_{n})^T$, where
\begin{equation*}
y_{j}=\left\{
\begin{array}{rl}
x_{ij}&~\text{if}~ v_{ij}\in V(G_i)\\
0&~\text{otherwise.}\\
\end{array}\right.
\end{equation*}
As $ e_{n_{i}}^{T}X=0 $ and coordinates of the  vector $Y$ corresponding to vertices in $\cup_{j\ne i}V_j$ of $H$ are zeros,  we have
\begin{align*}
A_{\alpha}(H)Y=\begin{pmatrix}
0\\
\vdots\\
0\\
\alpha (r_{i}+\alpha_i)X+(1-\alpha)\lambda_{ik} X\\
0\\
\vdots\\
0
\end{pmatrix}=\Big(\alpha(r_{i}+\alpha_i )+(1-\alpha)\lambda_{ik}\Big)Y.
\end{align*}
This shows that $Y$ is an eigenvector of $A_{\alpha}(H) $ corresponding to the eigenvalue $\alpha (r_{i}+\alpha_i)+(1-\alpha)\lambda_{ik}$, for every eigenvalue $\lambda_{ik}$, $2\leq k\leq n_i$, of $A(G_i)$. From this it follows that for $1\leq i\leq n$ and $2\leq k\leq n_i$, $\alpha (r_{i}+\alpha_i )+(1-\alpha)\lambda_{ik}$,  is an eigenvalue of $A_{\alpha}(H) $.  So, we have obtain $ \sum\limits_{i=1}^{n}n_{i}-n=N-n$  eigenvalues of $A_{\alpha}(H) $ with this procedure. The remaining $ n $ eigenvalues are the zeros of the characteristic polynomial of the following equitable quotient matrix.
\begin{equation*}
M =\begin{pmatrix}
\phi_{11}& (1-\alpha)\phi_{12}&\ldots & (1-\alpha)\phi_{1n}\\
(1-\alpha)\phi_{21}&\phi_{22}& \ldots & (1-\alpha)\phi_{2n}\\
\vdots &\vdots  &\ddots &\vdots\\
(1-\alpha)\phi_{n1}& (1-\alpha)\phi_{n2}&\ldots &\phi_{nn}
\end{pmatrix},
\end{equation*} 
where $\phi_{ii}=\alpha \alpha_i+r_i$, for $i=1,2,\dots,n$ and for $i\ne j$,  $\phi_{ij}=n_j$, if $v_i\sim v_j$, while as  $\phi_{ij}=0$, if $v_i\nsim v_j$.  This completes the proof.\qed

If in particular we take  $ G_{i}=K_{p_i} $ the complete graph on $p_i$ vertices (a ($p_i-1$)-regular graph), then using the fact that the adjacency eigenvalues of $K_{p_i}$ are $p_i-1$ with multiplicity $1$ and the eigenvalue $-1$ with multiplicity $p_i-1$, we have the following consequence of Theorem \ref{joined union}.
\begin{corollary}
Let $G$ be a graph of order $n\ge 3$. For $i=1,2,\dots,n$, let $ G_{i}=K_{p_{i}}$. Then the $ A_{\alpha} $-spectrum of $ G[K_{p_{1}},\ldots, K_{p_{n}}] $ consists of $ \alpha p_{i}-1+\alpha \alpha_i$, with multiplicity $ p_{i}-1 $, for each $ i=1,2,\dots,n $, where $\alpha_i=\sum\limits_{v_j\in N_{G}(v_{i})}n_{j}$. The remaining $n$ eigenvalues are given by the matrix 
\begin{equation*}
\begin{pmatrix}
\phi_{11}& (1-\alpha)\phi_{12}&\ldots & (1-\alpha)\phi_{1n}\\
(1-\alpha)\phi_{21}&\phi_{22}& \ldots & (1-\alpha)\phi_{2n}\\
\vdots &\vdots  &\ddots &\vdots\\
(1-\alpha)\phi_{n1}& (1-\alpha)\phi_{n2}&\ldots &\phi_{nn}
\end{pmatrix},
\end{equation*} 
where $\phi_{ii}=\alpha \alpha_i+p_i-1$, for $i=1,2,\dots,n$ and for $i\ne j$,  $\phi_{ij}=n_j$, if $v_i\sim v_j$, while as  $\phi_{ij}=0$, if $v_i\nsim v_j$.
\end{corollary}

We observe that the \textit{complete $p$-partite graph} $K_{n_1,n_2,\dots,n_p}$, is the joined union of graphs $ G_{i}=\overline{K}_{n_{i}} $, where underlying graph is $ K_p $. So, $K_{n_1,n_2,\dots,n_p}=K_{p}[\overline{K}_{n_1},\overline{K}_{n_2},\dots,\overline{K}_{n_p}].$
The next consequence of Theorem \ref{joined union} gives the $ A_{\alpha} $-spectrum of the complete $p$-partite graph $K_{n_1,n_2,\dots,n_p}$.

\begin{corollary}\label{p-partite}
The $ A_{\alpha} $-spectrum of the complete $p$-partite graph $K_{n_1,n_2,\dots,n_p}$ with $N=\displaystyle\sum_{i=1}^{p}n_{i}$ consists of the eigenvalue $\alpha(N-n_i)$, for $i=1,2,\dots,p$, each having multiplicity $n_i-1 $ and the remaining $ p $ eigenvalues are given by the matrix
\begin{equation*}
\begin{pmatrix}
\alpha(N-n_{1})& n_{2}(1-\alpha)&\ldots & n_{p}(1-\alpha)\\
n_{1}(1-\alpha)&\alpha(N-n_{2})& \ldots & n_{p}(1-\alpha)\\
\vdots &\vdots  &\ddots &\vdots\\
n_{1}(1-\alpha)& n_{2}(1-\alpha)&\ldots &\alpha(N-n_{p})
\end{pmatrix}.
\end{equation*}
\end{corollary}
{\textbf{Proof.}} This follows from Theorem \ref{joined union} by taking $ G=K_{p}$, $G_{i}=\overline{K}_{n_{i}}$, $a(v_{i},v_{j})=J_{n_i\times n_j}$, for all $i,j$, $\alpha_i=\sum\limits_{v_j\in N(v_{i})}n_{i}=N-n_{i}, r_{i}=0$ and $\lambda_{ik}=0$, $1\leq k\leq n_{i}-1 $ for $ 1\leq i\leq p. $ \qed

\begin{corollary}
Let $ G= K_{t,t,\dots,t} $ be a complete $ p$-partite graph with $ N=p t .$ The $ A_{\alpha} $-spectrum of $ G $ consists of the eigenvalue $\alpha t (p-1) $ with multiplicity $ tp-p $, the eigenvalue $t(\alpha p-1)$ with multiplicity $p-1$ and the eigenvalue $ t(p-1)$ with multiplicity $1$.
\end{corollary}
\textbf{Proof.} By Corollary \ref{p-partite}, we see that $ \alpha t(p-1) $ is an eigenvalue with multiplicity $ pt-p $ and other eigenvalues are given by
\begin{equation*}
\begin{pmatrix}
\alpha t(p-1)& t(1-\alpha)&\ldots & t(1-\alpha)\\
t(1-\alpha)&\alpha t(p-1)& \ldots & t(1-\alpha)\\
\vdots &\vdots  &\ddots &\vdots\\
t(1-\alpha)& t(1-\alpha)&\ldots &\alpha t(p-1)
\end{pmatrix}.
\end{equation*}
This matrix has symmetry of the type \eqref{pe1}, and so by taking $ X$ and $\beta$ as empty matrices, $B=[\alpha t(p-1)]$ and $ C=[t(1-\alpha)] $ in Lemma \ref{pe1}, we get the eigenvalue $ t(\alpha p-1) $ with multiplicity $ p-1 $. Finally $ B+(c-1)C =\alpha t(p-1)+(p-1)t(1-\alpha)=t(p-1)$ is the remaining eigenvalue. \qed

 The following observation gives the  $ A_{\alpha} $-spectrum of the join of a regular graph with the union of two regular graphs of distinct degrees. 
\begin{corollary}\label{spectrum of 3 graphs}
For $ i=1,2,3$, let $ G_{i} $ be $ r_{i} $ regular graphs of orders $ n_{i} $ having adjacency eigenvalues $ \lambda_{i1}=r_{i}\geq \lambda_{i2}\geq \dots\geq\lambda_{in_{i}} $.  Then the $ A_{\alpha} $-spectrum of $G= G_{1}\triangledown (G_{2}\cup G_{3})$ consists of the eigenvalues $\alpha(r_1+n-n_1)+(1-\alpha)\lambda_{1k}(A(G_{1}))$, $k=2,\dots,n_{1}$, the eigenvalues  $\alpha(n_1+r_2)+(1-\alpha)\lambda_{2k}(A(G_{2}))$, $k=2,\dots, n_{2}$, the eigenvalues  $\alpha(n_1+r_{3})+(1-\alpha)\lambda_{3k}(A(G_{3}))$, $k=2,\dots, n_{3}$, 
 where $n=n_1+n_2+n_3$.
 The remaining three eigenvalues are given by the  matrix
\begin{equation*}
\begin{pmatrix}
r_{2}+\alpha n_{1} & (1-\alpha)n_{1} &0\\
(1-\alpha)n_{2} & r_{1}+\alpha (n_{2}+n_{3}) &(1-\alpha)n_{3}\\
0 &(1-\alpha)n_{1} & r_{3}+\alpha n_{1}
\end{pmatrix}.
\end{equation*}
\end{corollary}
{\bf Proof.} Let $P_3$ be the path on $3$ vertices and let $G= G_{1}\triangledown (G_{2}\cup G_{3})$ be the join of the graphs $G_1$ and $G_2\cup G_3$. It is easy to see that $G=P_3[G_2,G_1,G_3]$, that is, $G$ is the joined union of the graphs $G_2, G_1, G_3$ when the parent graph is the path $P_3$. Now, using Theorem \ref{joined union} and noting that $\alpha_2=n_1$, $\alpha_1=n_2+n_3$ and $\alpha_3=n_1$, the result follows. \qed

\indent Let $G_1$ and $G_2$ be two graphs of order $n_1$ and $n_2$, respectively. Let $G=G_1\triangledown G_2$ be the join of the graphs $G_1$ and $G_2$. Clearly, $G=K_2[G_1,G_2]$, that is, $G$ is the joined union of $G_1$ and $G_2$ when the parent graph is $K_2$. The following consequence of Theorem  \ref{joined union}, gives the $ A_{\alpha} $-spectrum of the join of two regular graphs.
\begin{corollary}\label{join of two graphs}
Let $ G_{i}$ be an $ r_{i} $ regular graph of order $ n_{i} $ for $i=1,2$. Let $ \lambda_{ik}, 2\leq k\leq n_{i}, i=1,2 $ be the adjacency eigenvalues of $ G_{i} $. Then the $ A_{\alpha} $-spectrum of $ G=G_{1}\triangledown G_{2} $ consists of eigenvalues $ \alpha(r_{1}+n_{2})+(1-\alpha)\lambda_{1k}A(G_{1}) $, $k=2,\dots, n_{1} $, the eigenvalues $ \alpha(r_{2}+n_{1})+(1-\alpha)\lambda_{2k}A(G_{1}), k=2,\dots, n_{2} $ and the remaining two eigenvalues are given by the matrix
\begin{equation} \label{eq21}
\begin{pmatrix}
r_{1}+\alpha n_{2}& (1-\alpha)n_{2}\\
(1-\alpha)n_{1}& r_{2}+\alpha n_{1}
\end{pmatrix}.
\end{equation}
\end{corollary}

A \emph{friendship} graph $ F_{n} $ is a graph of order $ 2n+1 $, obtained by joining $ K_{1} $ with $ n $ copies of $ K_{2} $, that is, $ F_{n}=K_{1}\triangledown (nK_{2}) =K_{1,n}[K_{1},\underbrace{K_{2},K_{2},\dots,K_{2}}_{n}]$, where $K_1$ corresponds to the root vertex(vertex of degree greater than one) in $K_{1,n}$.  In particular, replacing some of $ K_{2} $'s by $ K_{1} $'s, in $F_n$ we get \emph{firefly} type graph denoted by $ F_{p,n-p}$ and written as 
\begin{equation*}
F_{p,n-p}=K_{1,n}[K_{1},\underbrace{K_{1},K_{1},\dots,K_{1}}_{p}\underbrace{K_{2},K_{2},\dots,K_{2}}_{n-p}].
\end{equation*}

The $ A_{\alpha} $-spectrum of the \emph{friendship} graph $ F_{n} $ and the \emph{firefly} type graph $ F_{p,n-p}$ are given by following corollary.
\begin{corollary}
\begin{itemize}
\item[\bf{(i)}] The $ A_{\alpha} $-spectrum of $ F_{n} $ consists of the eigenvalue $ 3\alpha -1 $ with multiplicity $ n$, the  remaining eigenvalues are the zeros of the  matrix
\begin{equation*}
\begin{pmatrix}
2\alpha n & 2(1-\alpha)& \dots & 2(1-\alpha)\\
1-\alpha & 1+\alpha & \dots & 2(1-\alpha)\\
\vdots & \vdots & \vdots & \vdots\\
1-\alpha & 2(1-\alpha) & \dots & 1+\alpha
\end{pmatrix}.
\end{equation*}
\item[ \bf{(ii)}] The $ A_{\alpha} $-spectrum of $ F_{p,n-p} $ consists of the eigenvalue $ 3\alpha-1 $ with multiplicity $ n-p $ and the eigenvalues of the matrix
\begin{equation*}
\begin{pmatrix}
\alpha (2n-p) & 1-\alpha& \dots & 2(1-\alpha)& 2(1-\alpha)& \dots & 2(1-\alpha)\\
1-\alpha & \alpha & \dots & 1-\alpha & 2(1-\alpha) &\dots & 2(1-\alpha)\\
\vdots & \vdots & \dots & \vdots & \vdots &  \dots & \vdots\\
1-\alpha & 1-\alpha & \dots & \alpha & 2(1-\alpha) & \dots & 2(1-\alpha)\\
1-\alpha & 1-\alpha & \dots & 1-\alpha & 1+\alpha & \dots & 2(1-\alpha)\\
\vdots & \vdots & \dots & \vdots & \vdots &  \dots & \vdots\\
1-\alpha & 1-\alpha & \dots & 1-\alpha & 2(1-\alpha) & \dots & 1+\alpha\\
\end{pmatrix}.
\end{equation*}
\end{itemize}
\end{corollary}

Let $ K_{a,b} $ be the complete bipartite graph with partite sets of cardinality $ a $ and $ b .$ Clearly $ K_{a,b}=\overline{K}_{a}\triangledown \overline{K}_{b} $. A \emph{complete split graph}, denoted by $CS_{\omega,n-\omega}$, is the graph consisting of a clique on $\omega$ vertices and an independent set (a subset of vertices of a graph is said to be an independent set if the subgraph induced by them is an empty graph) on the remaining $n-\omega$ vertices, such that each vertex of the clique is adjacent to every vertex of the independent set. Further, note that $CS_{\omega,n-\omega}=K_{\omega}\triangledown\overline{K}_{n-\omega
}$. Similarly, the \emph{cone graph} is given by $ C_{a,b}=C_{a}\triangledown \overline{K}_{b} $ and for $ b=1 $, we have \emph{wheel graph} $ W_{n+1}=C_{n}\triangledown K_{1} $ on $ n+1 $ vertices. As a special case of Corollary \ref{join of two graphs}, we compute the $ A_{\alpha} $-spectrum of the complete bipartite graph $ K_{a,b}$, the complete split graph $ CS_{\omega,n-\omega},$ the cone graph $ C_{a,b} $ and $ W_{n+1}$ in the following corollary.

\begin{corollary}\label{spec of join of two graphs}
\begin{itemize}
\item[\bf{(i)}] The $ A_{\alpha} $-spectrum of $ K_{a,b} $ consists of the eigenvalue $ \alpha b $ with multiplicity $ a-1$, the eigenvalue $\alpha a$ with multiplicity $ b-1 $ and the eigenvalues
\begin{equation*}
\dfrac{1}{2}\left ( \alpha(a+b)\pm \sqrt{\alpha^{2}(a+b)^2+4ab(1-2\alpha)} \right ).
\end{equation*}
\item[\bf{(ii)}] The $ A_{\alpha} $-spectrum of $CS_{\omega,n-\omega}$ is given by 
\begin{equation*}
\left\lbrace  (\alpha n-1)^{(\omega-1)}, (\alpha\omega)^{(n-\omega-1)}, \frac{1}{2}\left(5n-2\omega-6\pm \sqrt{(3(2\omega-n)-2(\omega-1))^{2}+4\omega(n-\omega)}\right)\right\rbrace .
\end{equation*}
\item[\bf{(iii)}] The $ A_{\alpha} $-spectrum of $ C_{a,b} $ consists of the eigenvalues $ \alpha(b+2)+2(1-\alpha)\cos\left( \frac{2\pi k}{b}\right) $, where $ k=2,\dots,b-1 $, the eigenvalue $ \alpha a $ with multiplicity $ b-1 $ and the eigenvalues
\begin{equation*}
\dfrac{1}{2}\left ( 2+\alpha (a+b)\pm \sqrt{\alpha^2(a+b)^{2}+\alpha(4a-4b-8ab)+4ab+4} \right ).
\end{equation*}
\item[\bf{(iv)}] The $ A_{\alpha} $-spectrum of $ W_{n+1} $ consists of the eigenvalues $ 3\alpha+2(1-\alpha)\cos\left( \frac{2\pi k}{n}\right) $, where $ k=2,\dots,n-1 $, and the eigenvalues
\begin{equation*}
\dfrac{1}{2}\left ( 2+\alpha(n+1) \pm \sqrt{\alpha^2(n+1)^{2}+\alpha(4-12n)+4n+4} \right ).
\end{equation*}

\end{itemize}
\end{corollary}
\textbf{Proof.} $ \textbf{(i)} $. This follows from Corollary \ref{spec of join of two graphs},  by taking $ n_{1}=a,n_{2}=b, r_{1}=r_{2}=0 $ and $ \lambda_{1k}=0,$ for $ k=2,\dots,a$ and $ \lambda_{2k} =0$ for each $ k=2,\dots,b$.\\
$ \textbf{(ii)} $. We recall that the adjacency spectrum of $ K_{\omega} $ is $ \{ \omega-1,-1^{(\omega-1)} \} $. Now the result follows from Corollary \ref{spec of join of two graphs} by taking $ n_{1}=\omega, n_{2}=n-\omega, r_{1}=\omega-1, r_{2}=0, \lambda_{1k}=-1,$ for $ k=2,\dots,\omega $ and $ \lambda_{2k}=0$ for $k=2,3,\dots,n-\omega .$ \\
$ \textbf{(iii)} $. We note that the adjacency spectrum (see \cite{cds}) of $ C_{n} $ is $ \left\lbrace 2\cos\left( \frac{2\pi k}{n}\right) : k=1,2,\dots,n\right\rbrace  $. Now result follows from Corollary \ref{spec of join of two graphs}, by taking $ n_{1}=a,n_{2}=b, r_{1}=2, r_{2}=0 $ and $ \lambda_{1k} =2\cos\left( \frac{2\pi i}{n}\right )$ for $ k=2,3,\dots,n. $  \\
$ \textbf{(iv)} $ Particular case of $ \textbf{(iii)} $ with $ b=1 $.\qed

\section{$A_{\alpha}$-Spectrum of power graphs of certain finite groups}

Let $ \mathcal{G} $ be a finite group of order $ n $ with identity element $ e $. Chakrabarty et al. \cite{sen} defined the undirected power graph $ \mathcal{P}(\mathcal{G}) $ of a group $ \mathcal{G} $ as an undirected graph with vertex set as $ \mathcal{G} $ and two vertices $ x,y\in \mathcal{G} $ are adjacent if and only if one is the positive power of other, that is, $ x^{i}=y $ or $ y^{j}=x $, for positive integers $i,j$ with  $ 2\leq i,j\leq n $. For some recent work on power graphs, we refer to \cite{cameron1, sen, tamiza, mehreen, survey} and the references therein. The adjacency spectrum, the Laplacian and the signless Laplacian spectrum of power graphs of finite cyclic and dihedral groups has been investigated in  \cite{banerjee,sriparna,asma, mehreen1, panda}.

An integer $ d $ dividing $ n $ is written as $ d|n $, and is called proper divisor if $ 1<d<n .$ Let $ d_{1}, d_{2},\dots,d_{r} $ be the distinct proper divisors of $ n. $ Let $\mathbb{G}_{n}$ be a simple graph with vertex set $ \{d_{i}: d_{i}|n, d_i\notin\{1,n\}, ~ 1 \leq i \leq r\} $ and edge set $\{ d_{i}d_{j}: d_{i}|d_{j}, ~ 1 \leq i< j\leq r\} $. If the \emph{canonical decomposition} of $n$ is $ n=p_{1}^{n_{1}}p_{2}^{n_{2}}\dots p_{r}^{n_{r}} $, where $ r,n_{1},n_{2},\dots,n_{r} $ are positive integers and $ p_{1},p_{2},\dots,p_{r} $ are distinct prime numbers, then from elementary number theory, the number of divisors of $ n $ are $\prod\limits_{i=1}^{r}(n_{i}+1) $, so the order of the graph  $\mathbb{G}_{n}$ is  $ |V( \mathbb{G}_{n})|=\prod\limits_{i=1}^{r}(n_{i}+1)-2. $ If $ n $ is neither a prime power nor product of two distinct primes, then $\mathbb{G}_{n}$ is a connected graph, see \cite{bilal}.\\

The following result shows that $ \alpha n-1 $ is always an $ A_\alpha $-eigenvalue of the power graph $ \mathcal{P}(\mathcal{G})$.

\begin{theorem}\label{mul of alpha n of p(G)}
Let $ \mathcal{G} $ be a finite group of order $ n\geq 3 $. Then $ \alpha n-1 $ is always an $ A_{\alpha} $-eigenvalue of the power graph $ \mathcal{P}(\mathcal{G}) $ with multiplicity at least $ b-1 $, where $ b $ is the number of vertices which are adjacent to every other vertex of $ \mathcal{P}(\mathcal{G}) $.
\end{theorem}
\textbf{Proof.} Let $ B(\mathcal{G}) $ be the subset of the vertex set $ \mathcal{G} $ of power graph  $ \mathcal{P}(\mathcal{G}) $ consisting of the identity $ e $ and those vertices which are adjacent to every other vertex of $ \mathcal{P}(\mathcal{G}) $. Assume that the cardinality of $ B(\mathcal{G}) $ is $ b $, so the induced subgraph of $ \mathcal{G} $ induced by the vertices in  $ B(\mathcal{G}) $ is a complete graph $ K_{b}$. Further, using the definition of subset $ B(\mathcal{G}) $, we see that  $ \mathcal{P}(\mathcal{G}) \cong K_{b}\triangledown  \mathcal{P}(\mathcal{G}\setminus B(\mathcal{G})) $. Applying Corollary \ref{join of two graphs}, we find that
\begin{equation*}
\alpha(r_{1}+n_{2})+(1-\alpha)\lambda_{1k}(G_{1})=\alpha(b-1+n-b)+(1-\alpha)(-1)=\alpha n-1
\end{equation*}
is an $ A_{\alpha} $-eigenvalue of $\mathcal{P}(\mathcal{G}) $  with multiplicity $ b-1 $. Since $ \alpha n-1 $ can also be an eigenvalue of $ \mathcal{P}(\mathcal{G}\setminus B(\mathcal{G})) $, therefore it follows that $ \alpha n-1 $ is an $ A_{\alpha} $-eigenvalue of $ \mathcal{P}(\mathcal{G}) $ with multiplicity at least $ b-1 $. \qed

Taking in particular $\mathcal{G}=\mathbb{Z}_{n}$, a finite cyclic group of order $n$, we have the following observation.

\begin{corollary}\label{mul of alpha n of Z_n}
Let $ \mathbb{Z}_{n} $ be a finite cyclic group of order $n\ge 3$. Then $ \alpha n-1 $ is an $ A_{\alpha} $-eigenvalue of $ \mathcal{P}(\mathbb{Z}_{n}) $ with multiplicity at least $\phi(n) $,  where $ \phi $ is Euler's totient function.
\end{corollary}
\textbf{Proof.} Since identity $ 0 $ and the invertible elements of the group $ \mathbb{Z}_{n} $ are adjacent to every other vertex of $\mathbb{Z}_{n}$ in the power graph $ \mathcal{P}(\mathbb{Z}_{n}) $, the result follows by Theorem \ref{mul of alpha n of p(G)} and the fact that there are $ \phi(n) $ invertible elements of $ \mathbb{Z}_{n} $.\qed

It will be an interesting problem to characterise the power graphs for which equality holds in Theorem \ref{mul of alpha n of p(G)} and Corollary \ref{mul of alpha n of Z_n}. Therefore, we leave the following problem. \\

\begin{problem}\label{Problem 1}
Characterize the  power graphs $\mathcal{P}(\mathcal{G})$, where $ \mathcal{G} $ is a finite group of order $n\ge 3$ having  $\alpha n-1$ an $ A_{\alpha} $-eigenvalue with multiplicity $ b-1 $, where $ b $ is the number of vertices which are adjacent to every other vertex of $ \mathcal{P}(\mathcal{G}) $. Also, characterize the  power graphs $\mathcal{P}(\mathbb{Z}_{n})$ having  $\alpha n-1$ an $ A_{\alpha} $-eigenvalue with multiplicity $\phi(n)$.

\end{problem}

From Theorem \ref{mul of alpha n of p(G)}, we see that if $ \mathcal{P}(\mathcal{G}\setminus B(\mathcal{G}) $ is known, then the $ A_{\alpha} $-spectrum of $ \mathcal{P}(\mathcal{G}) $ can be completely determined. So, it will be interesting to study the structure of $ \mathcal{P}(\mathcal{G}\setminus B(\mathcal{G}) $, and looking for graph parameters related to it.

The following result \cite{mehreen} shows that the power graph of a finite cyclic group $ \mathbb{Z}_{n} $ can be written as the joined union each of whose component is a clique.

\begin{theorem}\label{mehreen}
If $ \mathbb{Z}_{n} $ is a finite cyclic group of order $n\ge 3$, then the power graph $\mathcal{P}(\mathbb{Z}_{n})$ is given by  
\begin{equation*}
\mathcal{P}(\mathbb{Z}_{n}) =K_{\phi(n)+1}\triangledown \mathbb{G}_{n}[K_{\phi(d_{1})},K_{\phi(d_{2})},\dots,K_{\phi(d_{r})}],
\end{equation*} 
where $\mathbb{G}_{n}$ is the graph of order $r$ defined above. 
\end{theorem}

Using Theorems \ref{joined union} and its corollaries, we can compute the $ A_{\alpha} $-spectrum of $ \mathcal{P}(\mathbb{Z}_{n}) $ in terms of the adjacency spectrum of $ K_{\omega} $ and zeros of the characteristic polynomial of the auxiliary matrix. We recall \cite{sen} that $ \mathcal{P}(\mathcal{G}) $ is a complete graph if and only if  group $ \mathcal{G} $ is a cyclic group of order $ n=p^{z}$, $z\in \mathbb{N} $, where $ p $ is prime.

We form a new connected graph $ H=K_{1}\triangledown \mathbb{G}_{n} $ which is of diameter at most two if $ \mathbb{G}_{n} $ is not complete, otherwise its diameter is $ 1 $. Now, we compute the $ A_{\alpha} $-spectrum of the power graph of $ \mathbb{Z}_{n} $ by using Theorems \ref{joined union} and Theorem \ref{mehreen}.

\begin{theorem}\label{zn}
The $ A_{\alpha} $-spectrum of $ \mathcal{P}(\mathbb{Z}_n) $ is
\begin{equation*}
\left\lbrace (\alpha n-1)^{(\phi(n))}, (\alpha(\phi(d_{1})+\alpha_2)-1)^{(\phi(d_{1})-1)},\dots, (\alpha(\phi(d_{r})+\alpha_{r+1})-1)^{(\phi(d_{r})-1)} \right\rbrace
\end{equation*}
and the eigenvalues of the matrix
\begin{equation}\label{quotient matrix of Z_n}
M=\begin{pmatrix}
z_1 & (1-\alpha)\phi(d_{1})& (1-\alpha)\phi(d_{2})  & \dots & (1-\alpha)\phi(d_{n})\\
(1-\alpha)(\phi(n)+1) & z_2&  (1-\alpha)\psi_{23} & \dots & (1-\alpha ) \psi_{2r}\\
\vdots & \vdots &\vdots  & \cdots & \vdots\\
(1-\alpha)(\phi(n)+1) & (1-\alpha) \psi_{r2}& (1-\alpha)\psi_{r3} & \dots & z_{r+1}
\end{pmatrix},
\end{equation}
where $ z_1=\phi(n)+\alpha(n-1-\phi(n)) $, $ z_{i}=\alpha \alpha_{i}+\phi(d_{i})-1, \alpha_{i}=\sum\limits_{v_j\in N_{H}(v_{i})}n_{j}$, for $ i=2,\dots,r+1 $ and for $i\ne j$,  $\psi_{ij}=n_j$, if $v_i\sim v_j$, while as  $\psi_{ij}=0$, if $v_i\nsim v_j$.
\end{theorem}
\textbf{Proof.} Let $ \mathcal{G}\cong \mathbb{Z}_{n} $ be a finite cyclic group of order $ n $. Since by definition of power graph, the identity element $ 0 $ and the generators (which are $ \phi(n) $ in number) of the group $ \mathbb{Z}_{n} $, are adjacent to every other vertex of $ \mathcal{P}(\mathbb{Z}_{n}) $. So, by Theorem \ref{mehreen}, we have
\begin{equation*}
\mathcal{P}(\mathbb{Z}_{n}) =K_{\phi(n)+1}\triangledown \mathbb{G}_{n}[K_{\phi(d_{1})},K_{\phi(d_{2})},\dots,K_{\phi(d_{r})}]= H[K_{\phi(n)+1},K_{\phi(d_{1})},K_{\phi(d_{2})},\dots,K_{\phi(d_{r})}],
\end{equation*}
where $ H=K_{1}\triangledown \mathbb{G}_{n} $ is the graph with vertex set $ \{v_{1},\dots, v_{r+1}\} $. Taking $ G_{1}=K_{\phi(n)+1}$ and $ G_{i}=K_{\phi(d_{i-1})} $, for $ i=2,\dots,r+1 $, in  Theorem \ref{joined union} and using the fact that the adjacency spectrum of $K_{t}$ consists of the eigenvalue $t-1$ with multiplicity $1$ and the eigenvalue $-1$ with multiplicity $t-1$, it follows that \begin{align*}
 \alpha (r_{1}+\alpha_{1})+(1-\alpha)\lambda_{1k}A(G_{1}) &=\alpha(\phi(n)+\phi(d_{1})+\dots+\phi(d_{r}))-1+\alpha\\
 &= \alpha(\phi(n)+n-1-\phi(n))-1+\alpha=\alpha n-1
\end{align*} is an $ A_{\alpha} $ eigenvalues of $ \mathcal{P}(\mathbb{Z}_n) $ with multiplicity $ \phi(n)$. Note that we have used the fact that vertex $ v_{1} $ of graph $ H $ is adjacent to every other vertex of $ H $ and $ \alpha_1=\sum\limits_{ d|n, d\ne 1,n}\phi(d)=n-1-\phi(n)$, as $ \sum\limits_{ d|s}\phi(d) =s.$ Similarly, we can show that $ (\alpha(\phi(d_{1})+\alpha_2)-1),\dots, (\alpha(\phi(d_{r})+\alpha_{r+1})-1) $ are the $ A_{\alpha} $-eigenvalues of $ \mathcal{P}(\mathbb{Z}_n) $ with multiplicities $ \phi(d_{1})-1, \dots, \phi(d_{r})-1 $ respectively.  The remaining $ A_{\alpha} $-eigenvalues are the zeros of the characteristic polynomial of the equitable quotient matrix $ M $ given by \eqref{quotient matrix of Z_n}.\qed

Clearly, from Theorem \ref{zn} all the $ A_{\alpha} $-eigenvalues of the power graph $ \mathcal{P}(\mathbb{Z}_{n})$ are completely determined except the $r+1$-eigenvalues, which are the eigenvalues of   the matrix $M$ defined in Theorem \ref{zn}. Further, it is also clear that the matrix $M$  depends upon the structure of the graph $\mathbb{G}_{n} $, which is not known in general. However, if we give some particular values to $n$, then it may be possible to know the structure of graph $\mathbb{G}_{n}$ and hence about the matrix $M$. This information can be helpful to determine the $r+1$ remaining $ A_{\alpha} $-eigenvalues of the power graph $ \mathcal{P}(\mathbb{Z}_{n})$.

We have following consequences of Theorem \ref{zn}.

\begin{corollary}\label{pz}
If $ n=p^{z} $, where $ p $ is prime and $ z $ is a positive integer, then the $ A_{\alpha} $-spectrum of $ \mathcal{P}(\mathbb{Z}_n) $ is $ \{ n-1, (\alpha n-1)^{(n-1)} \} $.
\end{corollary}
\textbf{Proof.} If $ n=p^z $, where $ p $ is prime and $ z\in \mathbb{N} $, then as shown in \cite{sen}, $ \mathcal{P}(\mathbb{Z}_n)$ is isomorphic to the complete graph $ K_{n}$ and hence the result follows.\qed

\begin{corollary}\label{pq}
If $ n=pq $, where $p$ and $q$ ($ p<q $) are primes, then the $ A_{\alpha} $-spectrum of $ \mathcal{P}(\mathbb{Z}_n) $ consists of eigenvalues $ \{ (\alpha n-1)^{(\phi(n))}, (\alpha(\phi(n)+\phi(p)+1)-1)^{(\phi(p)-1)}, (\alpha(\phi(n)+\phi(q)+1)-1)^{(\phi(q)-1)}  \} $ and the three eigenvalues of the matrix given by \eqref{Qmat pq}.
\end{corollary}
\textbf{Proof.} Let $ n=pq $, where $p$ and $q$, $ p<q $, are distinct primes. Clearly, there are $ \phi(pq) $ generators in $\mathbb{Z}_{pq}$. These generators and identity element are adjacent to every other vertex of $\mathcal{P}(\mathbb{Z}_{n})$. So, by \cite{tamiza}, we have
\[  \mathcal{P}(\mathbb{Z}_n)= \left (K_{\phi(p)}\cup K_{\phi(q)}\right )\triangledown K_{\phi(pq)+1}=P_{3}[K_{\phi(p)},K_{\phi(pq)+1},K_{\phi(q)} ] . \]
Now, letting $ G_{1}=K_{\phi(pq)+1}, G_{2}=K_{\phi(p)} $ and $ G_{3}=K_{\phi(q)} $ in Corollary \ref{spectrum of 3 graphs}, we obtain the required $ A_{\alpha} $-eigenvalues and the remaining three eigenvalues are given by the matrix
\begin{equation} \label{Qmat pq}
\begin{pmatrix}
\phi(p)-1+\alpha(\phi(pq)+1) & (1-\alpha)(\phi(pq)+1) & 0\\
(1-\alpha)\phi(p) &\phi(pq)+\alpha(\phi(p)+\phi(q)) & (1-\alpha)\phi(q)\\
0 & (1-\alpha)(\phi(pq)+1) & \phi(q)-1+\alpha( \phi(pq)+1)
\end{pmatrix}.
\end{equation}\qed

Now, let $ n=pqr $, where $p,~ q,~ r$  with $ p<q<r$ are primes. From the definition of $\mathbb{G}_{n}$, the  vertex set and edge set of $\mathbb{G}_{n}$ are $\{ p, q, r, pq, pr, qr \}$ and $\{(p,pq),(p,pr),(q,pq),(q,qr),(r,pr),(r,qr)\} $. Let $ H=K_{1}\triangledown \mathbb{G}_{n}.$ Then $ \mathcal{P}(\mathbb{Z}_n)= H[K_{\phi(n)+1}, K_{\phi(d_{1})},\dots, K_{\phi(d_{6})}]$.

By Corollary \ref{mul of alpha n of Z_n}, $ \alpha n-1 $ is the $ A_{\alpha} $-eigenvalue with multiplicity $ \phi(n) $. By Theorem \ref{joined union}, other $ A_{\alpha} $-eigenvalues are $ \alpha(1+\phi(n)+\phi(p)+\phi(pr)+\phi(pq))-1,~ \alpha(1+\phi(n)+\phi(q)+\phi(qr)+\phi(pq))-1,~ \alpha(1+\phi(n)+\phi(r)+\phi(pr)+\phi(qr))-1,~ \alpha(1+\phi(n)+\phi(pq)+\phi(p)+\phi(q))-1,~ \alpha(1+\phi(n)+\phi(pr)+\phi(p)+\phi(r))-1, ~ \alpha(1+\phi(n)+\phi(qr)+\phi(q)+\phi(r))-1 $  with multiplicities $ \phi(p)-1,~ \phi(q)-1,~ \phi(r)-1,~ \phi(pq)-1,~ \phi(pr)-1,~ \phi(qr)-1$, respectively. The remaining $7$ eigenvalues are given by the following matrix
\begin{equation}\label{pqr}
\begin{pmatrix}
z_{1}& (1-\alpha)\phi(p) & (1-\alpha)\phi(q) & (1-\alpha)\phi(r)& (1-\alpha)\phi(pq)& (1-\alpha)\phi(pr)& (1-\alpha)\phi(qr)\\
c_1 & z_{2} & 0 & 0& (1-\alpha)\phi(pq) & (1-\alpha)\phi(pr)& 0 \\
c_1 & 0 & z_{3} & 0& (1-\alpha)\phi(pq) & 0 & (1-\alpha)\phi(qr) \\
c_1 & 0 & 0 & z_{4} & 0&  (1-\alpha)\phi(pr) & (1-\alpha)\phi(qr) \\
c_1 & (1-\alpha)\phi(p) & (1-\alpha)\phi(q) & 0 & z_{5} & 0  & 0  \\
c_1 & (1-\alpha)\phi(p) & 0 & (1-\alpha)\phi(r) & 0 & z_{6} & 0  \\
c_1 & 0 & (1-\alpha)\phi(q) & (1-\alpha)\phi(r) & 0 & 0 &  z_{7} 
\end{pmatrix},
\end{equation}
where, $ z_{1}=\phi(n)+\alpha(n-\phi(n)-1), ~ z_{2}= \phi(p)-1+\alpha(\phi(n)+1+\phi(pr)+\phi(pq)),~ z_{3}= \phi(q)-1+\alpha(\phi(n)+1+\phi(pq)+\phi(pr)), ~ z_{4}= \phi(r)-1+\alpha(\phi(n)+1+\phi(pr)+\phi(qr)), ~ z_{5}= \phi(pq)-1+\alpha(\phi(n)+1+\phi(p)+\phi(q)), ~ z_{6}= \phi(pr)-1+\alpha(\phi(n)+1+\phi(p)+\phi(r)), ~ z_{7}= \phi(qr)-1+\alpha(\phi(n)+1+\phi(q)+\phi(r))  $ and $ c_{1}=(1-\alpha)(\phi(n)+1) .$
\\

For the finite cyclic group $\mathbb{Z}_{pq^N}$, where $p$ and $q$ are distinct primes and $N$ is a positive integer, we have the following observation.

\begin{theorem} 
Let $ \mathcal{P}(\mathbb{Z}_{pq^N}) $ be the power graph of the finite cyclic group $\mathbb{Z}_{pq^N}$ of order $ n=pq^N $, where $p$ and $q$ are distinct primes and $N$ a positive integer. Then the following hold.
\begin{itemize}
\item[\bf{(i)}] If $ N $ is even, say $N=2m$. Then $ A_{\alpha} $-spectrum of $ \mathcal{P}(\mathbb{Z}_{pq^N})$ is 
\begin{align*}
\{&(\alpha n-1)^{(\phi(n))},\left (\alpha(\phi(n)+\phi(p)q^{2m-1}+1)-1\right )^{(\phi(p)-1)},(\alpha(\phi(n)+q^{2m}+\phi(p)(q^{2m-1}-1))\\
&-1)^{(\phi(q)-1)},  \dots,\left (\alpha(\phi(n)+q^{2m}+\phi(p)(q^{2m-1}-q^{m-1}))-1\right )^{(\phi(q^m)-1)},\dots, (\alpha(\phi(n)+q^{2m})\\
&-1)^{(\phi(q^{2m})-1)},\left (\alpha(\phi(n)+\phi(q)+\phi(p)q^{2m-1}+1)-1\right )^{(\phi(pq)-1)},\dots, (\alpha(\phi(n)+q^{m}+\phi(p)q^{2m-1})\\
&-1)^{(\phi(pq^m)-1)},\dots, \left (\alpha(\phi(n)+q^{2m-1}+\phi(p)q^{2m-1})-1\right )^{(\phi(pq^{2m-1})-1)} \} \end{align*} 
together with the zeros of the characteristic polynomial of the matrix given by \eqref{mat 1}.
\item[\bf{(ii)}] If $ N $ is odd, say $N=2m+1$. Then $ A_{\alpha} $-spectrum of $ \mathcal{P}(\mathbb{Z}_{pq^N})$ is 
\begin{align*}
\{(\alpha n&-1)^{(\phi(n))},\left (\alpha(\phi(n)+\phi(p)q^{2m}+1)-1\right )^{(\phi(p)-1)},\left (\alpha(\phi(n)+q^{2m+1}+\phi(p)(q^{2m}-1))-1\right  )^{(\phi(q)-1)},\\
\dots,&\left (\alpha(\phi(n)+q^{2m+1}+\phi(p)(q^{2m}-q^{m}))-1\right )^{(\phi(q^{m+1})-1)},\dots,\left (\alpha(\phi(n)+q^{2m+1})-1\right )^{(\phi(q^{2m+1})-1)},\\
&\left (\alpha(\phi(n)+\phi(q)+ \phi(p)q^{2m}+1)-1\right )^{(\phi(pq)-1)},\dots,\left (\alpha(\phi(n)+q^{m+1}+\phi(p)q^{2m})-1\right )^{(\phi(pq^{m+1})-1)},\\
\dots,& \left (\alpha(\phi(n)+q^{2m}+\phi(p)q^{2m-1})-1\right )^{(\phi(pq^{2m}-1)} \} \end{align*} 
together with the zeros of the characteristic polynomial of the matrix given by \eqref{mat 2}.\
\end{itemize}
\end{theorem}
\textbf{Proof.} $ \textbf{(i)} $. Let $N$ be even, say $ N=2m, m\in \mathbb{N} $. Clearly, the order of $ \mathbb{G}_{pq^{2m}} $ is $n= 4m$ and its vertex set is $ \{p, q,q^2,\dots,q^m,\dots,q^{2m}, pq, pq^2,\dots,pq^m,\dots,pq^{2m-1}\} $, where $ p $ and $ q $ are distinct primes. Also, in the graph $ \mathbb{G}_{pq^N} $, we have $ p\sim pq^i $, for $ i=1,2,\dots, 2m-1 $, $ q^i\sim q^{j},$ for all $1\le i<j\le 2m$, $q^i\sim pq^{j} $, for all  $ 1\le i\leq j\le 2m-1$ and  $ pq^i\sim pq^{j}$ for all $ 1\le i<j\le 2m-1$, where we are avoiding the divisor which are outside the vertex set of $ \mathbb{G}_{pq^{2m}}$. Thus, by Theorem \ref{mehreen}, we have \begin{align*}
\mathcal{P}(\mathbb{Z}_{pq^{N}})&= K_{\phi(pq^{N})+1}\triangledown \mathbb{G}_{pq^{N}}[K_{\phi(p)},K_{\phi(q)},\dots,K_{\phi(q^m)},\dots,K_{\phi(q^{2m})},K_{\phi(pq)},\dots,K_{\phi(pq^m)},\dots,K_{\phi(pq^{2m-1})}]\\
&=H[K_{\phi(pq^{N})+1},K_{\phi(p)},K_{\phi(q)},\dots,K_{\phi(q^m)},\dots,K_{\phi(q^{2m})},K_{\phi(pq)},\dots,K_{\phi(pq^m)},\dots,K_{\phi(pq^{2m-1})}],\end{align*}
where $ H=K_{1}\triangledown \mathbb{G}_{pq^N} $. By Corollary \ref{mul of alpha n of Z_n}, we find that $ \alpha n-1 $ is an $ A_{\alpha} $-eigenvalue with multiplicity $ \phi(n) $. Now, by using the Theorem \ref{joined union} with $ n_{1}=\phi(pq^{2m})+1, ~ n_{2}=\phi(p), ~ n_{i}=\phi(q^j), i=3,\dots,2m+2, ~ j=1,2,\dots, 2m, ~ n_{i}=\phi(pq^j), i=2m+3,\dots,4m+1, ~ j=1,2,\dots,2m-1 $, it follows that \begin{align*}
 \alpha(\phi(p)&-1+\phi(n)+1+\phi(pq)+\dots+\phi(pq^m)+\dots+\phi(pq^{2m-1}))+(1-\alpha)(-1)\\
&= \alpha(\phi(n)+\phi(p)(1+\phi(q)+\dots+\phi(q^{2m-1})))-1+\alpha=\alpha(\phi(n)+\phi(p)q^{2m-1}+1)-1  
\end{align*} is the $ A_{\alpha} $-eigenvalue with multiplicity $ \phi(p)-1 $. Note that we have used the fact that $\phi(xy)=\phi(x)\phi(y),$ provided that g.c.d. of $x,y$ is $1$ and $\sum\limits_{i=1}^{k}\phi(p^i)=p^k-1$. Using the similar procedure we get the other $ A_{\alpha} $-eigenvalues as given in the statement. The remaining $ A_{\alpha} $-eigenvalues are given by the quotient matrix
\begin{equation}\label{mat 1}
\left( {\begin{array}{c c c c c c c c c c c c c c}
z_{1} & n_{2}^{'}& n_{3}^{'} & \cdots & n_{m+2}^{'} & n_{m+3}^{'} & \cdots & n_{2m+2}^{'} & n_{2m+3}^{'} & \cdots & n_{3m+2}^{'} &n_{3m+3}^{'} & \cdots & n_{4m+1}^{'}\\
n_{1}^{'} & z_{2}& 0 & \cdots & 0 & 0 & \cdots & 0 & n_{2m+3}^{'} & \cdots & n_{3m+2}^{'} &n_{3m+3}^{'} & \cdots & n_{4m+1}^{'}\\
n_{1}^{'} & 0& z_{3} & \cdots & n_{m+2}^{'} & n_{m+3}^{'} & \cdots & n_{2m+2}^{'} & n_{2m+3}^{'} & \cdots & n_{3m+2}^{'} &n_{3m+3}^{'} & \cdots & n_{4m+1}^{'}\\
\vdots & \vdots & \vdots & \vdots & \vdots & \vdots & \vdots & \vdots & \vdots & \vdots & \vdots & \vdots & \vdots & \vdots \\
n_{1}^{'} & 0& n_{3}^{'} & \cdots & z_{m+2} & n_{m+3}^{'} & \cdots & n_{2m+2}^{'} & 0 & \cdots & n_{3m+2}^{'} &n_{3m+3}^{'} & \cdots & n_{4m+1}^{'}\\
n_{1}^{'} & n_{2}^{'}& n_{3}^{'} & \cdots & n_{m+2}^{'} & z_{m+3} & \cdots & n_{2m+2}^{'} & 0 & \cdots & 0 & 0 & \cdots & 0 \\
\vdots & \vdots & \vdots & \vdots & \vdots & \vdots & \vdots & \vdots & \vdots & \vdots & \vdots & \vdots & \vdots & \vdots \\
n_{1}^{'} & n_{2}^{'}& n_{3}^{'} & \cdots & n_{m+2}^{'} & n_{m+3}^{'} & \cdots & z_{2m+2}^{'} & 0 & \cdots & 0 & 0 & \cdots & 0\\
n_{1}^{'} & n_{2}^{'}& n_{3}^{'} & \cdots & n_{m+2}^{'} & n_{m+3}^{'} & \cdots & n_{2m+2}^{'} & z_{2m+3}^{'} & \cdots & n_{3m+2}^{'} &n_{3m+3}^{'} & \cdots & n_{4m+1}^{'}\\
\vdots &\vdots & \vdots & \vdots & \vdots & \vdots & \vdots & \vdots & \vdots & \vdots & \vdots & \vdots & \vdots & \vdots \\
n_{1}^{'} & n_{2}^{'}& n_{3}^{'} & \cdots & 0 & 0 & \cdots & 0 & n_{2m+3}^{'} & \cdots & z_{3m+2}^{'} &n_{3m+3}^{'} & \cdots & n_{4m+1}^{'}\\
n_{1}^{'} & n_{2}^{'}& n_{3}^{'} & \cdots & n_{m+2}^{'} & n_{m+3}^{'} & \cdots & 0 & n_{2m+3}^{'} & \cdots & n_{3m+2}^{'} &z_{3m+3}^{'} & \cdots & n_{4m+1}^{'}\\
\vdots & \vdots & \vdots & \vdots & \vdots & \vdots & \vdots & \vdots & \vdots & \vdots & \vdots & \vdots & \vdots & \vdots \\
n_{1}^{'} & n_{2}^{'}& n_{3}^{'} & \cdots & n_{m+2}^{'} & n_{m+3}^{'} & \cdots & 0 & n_{2m+3}^{'} & \cdots & n_{3m+2}^{'} &n_{3m+3}^{'} & \cdots & z_{4m+1}^{'}\\
\end{array}} \right),
\end{equation}
where $ n_{i}^{'}=(1-\alpha)n_{i},~ \text{for} ~ i=1,2,\dots,4m+1 $ and $ z_{1}=\phi(n)+\alpha(n-\phi(n)-1), ~ z_{2}=\phi(p)-1+\alpha(\phi(n)+1+\phi(p)(q^{2m-1}-q)), ~ z_{3}=\phi(q)-1+\alpha(\phi(n)+1+ q^{2m}-q+\phi(p)(q^{2m-1}-1)) , \dots ,z_{m+2}=\phi(q^{m})-1+\alpha(\phi(n)+q^{2m}-\phi(q^m) +\phi(p)(q^{2m}-q^{m-1})), z_{m+3}=\phi(q^{m+1})-1+\alpha(\phi(n)+q^{2m}-\phi(q^{m+1}) +\phi(p)(q^{2m}-q^{m})), \dots, z_{2m+2}=\phi(q^{2m})-1+\alpha(\phi(n)+q^{2m}-\phi(q^{2m})),z_{2m+3}=\phi(pq)-1+\alpha(\phi(n)+\phi(p)+\phi(q)+\phi(p)(q^{2m-1}-q)), \dots, z_{3m+2}=\phi(pq^{m})-1+\alpha(\phi(n)+\phi(p)+q^{m}+\phi(p)(q^{2m-1}-1-\phi(q^{m}))), z_{3m+3}=\phi(pq^{m+1})-1+\alpha(\phi(n)+\phi(p)+q^{m+1}+\phi(p)(q^{2m}-1-\phi(q^{m+1}))),\dots, z_{4m+1}=\phi(pq^{2m-1})-1+\alpha(\phi(n)+q^{2m-1}+\phi(p)q^{2m-2})  .$

\noindent$ \textbf{(ii)} $. Now, let $N$ be odd, say $ N=2m+1, m\in \mathbb{N} $. Clearly, the order of $ \mathbb{G}_{pq^{2m+1}} $  is $ 4m+2$ and its vertex set is $ \{p, q,q^2,\dots,q^{m+1},\dots,q^{2m+1}, pq, pq^2,\dots,pq^{m+1},\dots,pq^{2m}\} $, where $ p $ and $ q $ are distinct primes. Also, in $ \mathbb{G}_{pq^{2m+1}} $, we have $ p\sim pq^i $, for $ i=1,2,\dots, N-1 $, $ q^i\sim q^{j}$, for all $1\le i<j\le 2m+1$, $q^{i}\sim pq^{j} $, for all $1\le i\le j\le 2m$ and $ pq^i\sim pq^{j},$ for all $1\le i<j\le 2m$. Note that we are avoiding the divisor which are outside the vertex set of $ \mathbb{G}_{pq^{2m+1}} $. Therefore, by Theorem \ref{mehreen}, we have 
\begin{align*}
 \mathcal{P}(\mathbb{Z}_{pq^{N}})= K_{\phi(pq^{N})+1}&\triangledown \mathbb{G}_{pq^{N}}[K_{\phi(p)},K_{\phi(q)},K_{\phi(q^2)},\dots,K_{\phi(q^{m+1})},\dots,K_{\phi(q^{N})},K_{\phi(pq)},\dots,K_{\phi(pq^{m+1})},\dots,\\
&K_{\phi(pq^{N-1})}]=H[K_{\phi(pq^{N})+1},K_{\phi(p)},K_{\phi(q)},K_{\phi(q^2)},\dots,K_{\phi(q^{m+1})},\dots,K_{\phi(q^{N})},K_{\phi(pq)},\\
&\dots,K_{\phi(pq^{m+1})},\dots,K_{\phi(pq^{N-1})}] ,\end{align*}
where $ H=K_{1}\triangledown \mathbb{G}_{pq^N} $.  By Corollary \ref{mul of alpha n of Z_n}, it follows that $ \alpha n-1 $ is an $ A_{\alpha} $-eigenvalue with multiplicity $ \phi(n) $. Using  Theorem \ref{joined union} with $ n_{1}=\phi(pq^{2m+1})+1, ~ n_{2}=\phi(p), ~ n_{i}=\phi(q^j), i=2,\dots,2m+3, ~ j=1,2,\dots, 2m+1, ~ n_{i}=\phi(pq^j), i=2m+3,\dots,4m+2, ~ j=1,2,\dots,2m $, it follows that $\alpha(\phi(n)+\phi(p)+\phi(p)(q^{2m}-1)+1)-1  $ is the $ A_{\alpha} $-eigenvalue with multiplicity $ \phi(p)-1 $. Proceeding in a similar way, we get the other $ A_{\alpha} $-eigenvalues as given in statement. The remaining $ A_{\alpha} $-eigenvalues are given by the following quotient matrix
\begin{equation}\label{mat 2}
\left( {\begin{array}{c c c c c c c c c c c c c c}
z_{1} & n_{2}^{'}& n_{3}^{'} & \cdots & n_{m+3}^{'} & n_{m+4}^{'} & \cdots & n_{2m+3}^{'} & n_{2m+4}^{'} & \cdots & n_{3m+3}^{'} &n_{3m+4}^{'} & \cdots & n_{4m+2}^{'}\\
n_{1}^{'} & z_{2}& 0 & \cdots & 0 & 0 & \cdots & 0 & n_{2m+4}^{'} & \cdots & n_{3m+3}^{'} &n_{3m+4}^{'} & \cdots & n_{4m+2}^{'}\\
n_{1}^{'} & 0& z_{3} & \cdots & n_{m+3}^{'} & n_{m+4}^{'} & \cdots & n_{2m+3}^{'} & n_{2m+4}^{'} & \cdots & n_{3m+3}^{'} &n_{3m+4}^{'} & \cdots & n_{4m+2}^{'}\\
\vdots & \vdots & \vdots & \vdots & \vdots & \vdots & \vdots & \vdots & \vdots & \vdots & \vdots & \vdots & \vdots & \vdots \\
n_{1}^{'} & 0& n_{3}^{'} & \cdots & z_{m+3} & n_{m+4}^{'} & \cdots & n_{2m+3}^{'} & 0 & \cdots & n_{3m+3}^{'} &n_{3m+4}^{'} & \cdots & n_{4m+2}^{'}\\
n_{1}^{'} & n_{2}^{'}& n_{3}^{'} & \cdots & n_{m+3}^{'} & z_{m+4} & \cdots & n_{2m+3}^{'} & 0 & \cdots & 0 & 0 & \cdots & 0 \\
\vdots & \vdots & \vdots & \vdots & \vdots & \vdots & \vdots & \vdots & \vdots & \vdots & \vdots & \vdots & \vdots & \vdots \\
n_{1}^{'} & n_{2}^{'}& n_{3}^{'} & \cdots & n_{m+3}^{'} & n_{m+4}^{'} & \cdots & z_{2m+3}^{'} & 0 & \cdots & 0 & 0 & \cdots & 0\\
n_{1}^{'} & n_{2}^{'}& n_{3}^{'} & \cdots & n_{m+3}^{'} & n_{m+4}^{'} & \cdots & n_{2m+3}^{'} & z_{2m+4}^{'} & \cdots & n_{3m+3}^{'} &n_{3m+4}^{'} & \cdots & n_{4m+2}^{'}\\
\vdots &\vdots & \vdots & \vdots & \vdots & \vdots & \vdots & \vdots & \vdots & \vdots & \vdots & \vdots & \vdots & \vdots \\
n_{1}^{'} & n_{2}^{'}& n_{3}^{'} & \cdots & 0 & 0 & \cdots & 0 & n_{2m+4}^{'} & \cdots & z_{3m+3}^{'} &n_{3m+4}^{'} & \cdots & n_{4m+2}^{'}\\
n_{1}^{'} & n_{2}^{'}& n_{3}^{'} & \cdots & n_{m+3}^{'} & n_{m+4}^{'} & \cdots & 0 & n_{2m+4}^{'} & \cdots & n_{3m+3}^{'} &z_{3m+4}^{'} & \cdots & n_{4m+2}^{'}\\
\vdots & \vdots & \vdots & \vdots & \vdots & \vdots & \vdots & \vdots & \vdots & \vdots & \vdots & \vdots & \vdots & \vdots \\
n_{1}^{'} & n_{2}^{'}& n_{3}^{'} & \cdots & n_{m+3}^{'} & n_{m+4}^{'} & \cdots & 0 & n_{2m+4}^{'} & \cdots & n_{3m+3}^{'} &n_{3m+4}^{'} & \cdots & z_{4m+2}^{'}\\
\end{array}} \right),
\end{equation}
where $ n_{i}^{'}=(1-\alpha)n_{i},~ \text{for} ~ i=1,2,\dots,4m+1 $ and $ z_{1}=\phi(n)+\alpha(n-\phi(n)-1), ~ z_{2}=\phi(p)-1+\alpha(\phi(n)+1+\phi(p)(q^{2m}-1)), ~ z_{3}=\phi(q)-1+\alpha(\phi(n)+1+ q^{2m+1}-q+\phi(p)(q^{2m}-1)) , \dots ,z_{m+3}=\phi(q^{m+1})-1+\alpha(\phi(n)+q^{2m+1}-\phi(q^{m+1}) +\phi(p)(q^{2m}-q^{m})), z_{m+4}=\phi(q^{m+2})-1+\alpha(\phi(n)+q^{2m+1}-\phi(q^{m+2}) +\phi(p)(q^{2m}-q^{m+1})), \dots, z_{2m+3}=\phi(q^{2m+1})-1+\alpha(\phi(n)+q^{2m+1}),z_{2m+4}=\phi(pq)-1+\alpha(\phi(n)+\phi(p)+\phi(q)+\phi(p)(q^{2m}-q)), \dots, z_{3m+3}=\phi(pq^{m+1})-1+\alpha(\phi(n)+\phi(p)+q^{m+1}+\phi(p)(q^{2m}-1-\phi(q^{m+1}))), z_{3m+4}=\phi(pq^{m+2})-1+\alpha(\phi(n)+\phi(p)+q^{m+2}+\phi(p)(q^{2m}-1-\phi(q^{m+2}))),\dots, z_{4m+2}=\phi(pq^{2m})-1+\alpha(\phi(n)+q^{2m}+\phi(p)q^{2m-1}).$ This completes the proof.\qed

Although for $n=p^{z}, pq,pqr, pq^{N}$, where $p,q,r$ are distinct primes and $z, N$ are positive integers, we have described the $A_{\alpha}$-spectrum of the power graph $\mathcal{P}(\mathbb{Z}_{n})$, the following problem will be of interest.

\begin{problem}
 Find the $ A_{\alpha} $-spectrum of the power graph $\mathcal{P}(\mathbb{Z}_{n}) $ of the finite cyclic group $\mathbb{Z}_{n}$, when $ n\in\{(pq)^2,~pqrs,~(pq)^m,~ m\in \mathbb{N}\} $ and generalize (if possible) for any $ n $.
\end{problem}

A group $ \mathcal{G} $ is an elementary abelian $ p$-group, if each non trivial element has order $ p $, where $ p $ is prime. In order to find $ A_{\alpha} $-spectrum of finite elementary abelian groups of prime power order, we need following lemma \cite{tamiza}.
\begin{lemma}\label{pgroup} Let $ \mathcal{G} $ be an elementary abelian group  of order $ p^n $ for some prime number $ p $ and positive integer $ n .$ Then $ \mathcal{P}(\mathcal{G})\cong K_{1}\triangledown \left ( \cup_{i=1}^{l} K_{p-1} \right ), $ where $ l=\frac{p^n-1}{p-1}. $
\end{lemma}

The following result gives the $ A_{\alpha} $-spectrum of power graph of elementary $ p$-group of order $ p^{n} $.
\begin{theorem}\label{pgroupt}
Let $ \mathcal{G} $ be an elementary abelian $ p$-group of order $ p^n $. The $ A_{\alpha} $-spectrum of $ \mathcal{P}(\mathcal{G}) $ consists of the eigenvalue $ \alpha p-1 $ with multiplicity $ l(p-2) $, where $ l=\frac{p^n-1}{p-1}$ and the eigenvalues of the matrix given by \eqref{Qmat of elementary abelian group}.
\end{theorem}
\textbf{Proof.} Let $ \mathcal{G} $ be an elementary abelian $ p$-group of order $ p^n $ for some prime $ p $ and positive integer $ n .$ Then, by Lemma \ref{pgroup}, we have  $ \mathcal{P}(\mathcal{G}) \cong K_{1}\triangledown \left (\cup_{i=1}^{l} K_{p-1} \right ) $ and so it follows that
\[\mathcal{P}(\mathcal{G}) \cong S[K_{1},\underbrace{ K_{p-1},\dots,K_{p-1}}_{l}], \] 
where $ K_{p-1} $ occurs $ l=\frac{p^n-1}{p-1} $ times and $ S=K_{1}\triangledown \overline{K}_{l}=K_{1,l}. $ Since $ G_{1}=K_{1}, ~ G_{i}=K_{p-1} $, for $ i=2,\dots,l+1 $ and $ \alpha_1=l(p-1)=p^n-1 $, $ \alpha_2=\dots=\alpha_{l+1}= 1$, therefore, using Theorem \ref{joined union}, we see that $ \alpha(r_{2}+\alpha_{2})+(1-\alpha)\lambda_{2k}A(G_{2})=\alpha(p-2+1)+\alpha-1=\alpha p-1 $ is the $ A_{\alpha} $-eigenvalue of $ \mathcal{P}(\mathcal{G}) $  with multiplicity $ l(p-2) $. The $l+1$ remaining $ A_{\alpha} $-eigenvalues of $ \mathcal{P}(\mathcal{G}) $ are given by the matrix 
\begin{equation}\label{Qmat of elementary abelian group}
\begin{pmatrix}
\alpha l(p-1) & (1-\alpha)(p-1)& (1-\alpha)(p-1) & \dots & (1-\alpha)(p-1)\\
1-\alpha & \alpha +p-2 & 0 & \dots & 0\\
\vdots & \vdots & \vdots & \cdots & \vdots\\
1-\alpha & 0 & 0 & \cdots & \alpha +p-2
\end{pmatrix}. 
\end{equation}
It is easy to see that $ \alpha+p-2 $ with multiplicity $ l-1 $ is the eigenvalue of the matrix given by \eqref{Qmat of elementary abelian group} and the other two eigenvalues are given by the following matrix 
\begin{equation*}
\begin{pmatrix}
\alpha l(p-1) & l (1-\alpha)(p-1)\\
1-\alpha & \alpha+p-2
\end{pmatrix},
\end{equation*} 
which are calculated as $ \alpha(lp+1-l)+p-2\pm\sqrt{(\alpha lp+p+\alpha-\alpha l-2)^2-4l(\alpha p^2-\alpha p-p+1)} $. This completes the proof.\qed

The following observation can be seen in \cite{tamiza}.
\begin{lemma}\label{pq non abelian}
Let $ \mathcal{G} $ be a finite group of order $ pq  $, where $ p $ and $ q $ are primes. Then $ \mathcal{G} $ is non abelian if and only if $ \mathcal{P}(\mathcal{G}) \cong K_{1}\triangledown (qK_{p-1}\cup K_{q-1})$.
\end{lemma}

Now, we obtain the $ A_{\alpha} $-spectrum of the non abelian group whose order is the product of two distinct primes. 

\begin{theorem}\label{pq non abelian me}
Let $ \mathcal{G} $ be a non abelian group of order $ n =pq$, where $p,~q$ ($ p<q $) are primes. Then the $ A_{\alpha} $-spectrum of $ \mathcal{P}(\mathcal{G}) $ consists of the eigenvalue $ \alpha p-1 $ with multiplicity $ q(p-2) $, the eigenvalue $ \alpha q-1 $ with multiplicity $ q-2 $ and the eigenvalues of the matrix given by \eqref{Qmat of pq non abelian group}.
\end{theorem}
\textbf{Proof.} Let $ \mathcal{G} $ be a non abelian group of order $ pq $. By Lemma \ref{pq non abelian}, the power graph of $  \mathcal{P}(\mathcal{G}) $ can be written as \[ \mathcal{P}(\mathcal{G})=K_{1,q+1}[K_{1},\underbrace{K_{p-1},K_{p-1},\dots,K_{p-1}}_{q},K_{q-1}]. \]  Since $ \alpha_1= pq-1,~ \alpha_2=\dots=\alpha_{q+2}= 1 $, so by Theorem \ref{joined union}, $ \alpha p-1 $ is the $ A_{\alpha} $-eigenvalue with multiplicity $ q(p-1) $. Similarly $\alpha q-1 $ is the $ A_{\alpha} $-eigenvalue with multiplicity $ q-2 $ and the remaining $ q+2 $ $ A_{\alpha} $-eigenvalues of $ \mathcal{P}(\mathcal{G}) $ are the eigenvalues of the matrix $ M $ given in \eqref{Qmat of pq non abelian group}.
\begin{equation}\label{Qmat of pq non abelian group}
M=\begin{pmatrix}
\alpha (pq-1) & (1-\alpha)(p-1)& (1-\alpha)(p-1) & \dots & (1-\alpha)(p-1) & (1-\alpha)(q-1)\\
1-\alpha & \alpha +p-2 & 0 & \dots & 0 & 0\\
\vdots & \vdots & \vdots & \cdots & \vdots & \vdots\\
1-\alpha & 0 & 0 & \cdots & \alpha +p-2 &0\\
1-\alpha & 0 & 0 & \cdots & 0 &\alpha+ q-2\\
\end{pmatrix}
\end{equation}
Clearly, $ \alpha+p-2 $ is an $ A_{\alpha} $-eigenvalue of matrix $ M $  with multiplicity $ q-1 $. The remaining three eigenvalues of matrix $ M $ are the eigenvalues of the following matrix
\begin{equation*}
\begin{pmatrix}
\alpha l(p-1) & q (1-\alpha)(p-1) & (1-\alpha)(q-1)\\
1-\alpha & \alpha+p-2 & 0\\
1-\alpha & 0 & \alpha+q-2
\end{pmatrix}.
\end{equation*}\qed

Next, we find the $ A_{\alpha} $-spectrum of the dihedral group and dicyclic group for some particular values of $ n $. The dihedral group of order $ 2n $ and the dicyclic group of order $ 4n $ are denoted  as
\begin{align*}
D_{2n}&=< a,b ~|~ a^{n}=b^{2}=e, bab=a^{-1}>,\\
Q_{n}&=<a,b ~|~ a^{2n}=e, b^2=a^n, ab=ba^{-1}>. 
\end{align*}
If $ n $ is a power of $ 2 ,$ then $ Q_{n} $ is called the \emph{generalized quaternion group} of order $ 4n. $

\begin{theorem}\label{dihedral}
If $ n $ is a prime power, then the $ A_{\alpha} $-spectrum of $ \mathcal{P}(D_{2n}) $ consists of the eigenvalue $ \alpha n-1 $ with multiplicity $ n-2 $, the eigenvalue $ \alpha $ with multiplicity $ n-1 $ and the zeros of following polynomial 
\[(\alpha -x) \left(\alpha ^2-2 \alpha +2 \alpha  n^2+\alpha ^2 n-3 \alpha  n-2 \alpha  n x-n x-n+x^2+2 x+2\right)-(1-\alpha )^2 n (\alpha +n-x-2). \]
\end{theorem}
\textbf{Proof.} Consider the dihedral group $D_{2n}=< a,b ~|~ a^{n}=b^{2}=e, bab=a^{-1}>$. Clearly, $ <a> $ generates a cyclic subgroup of order $ n $, which is isomorphic to $ \mathbb{Z}_{n} $. The remaining $ n $ elements of $ D_{2n} $ form an independent set of $ \mathcal{P}(D_{2n}) $, sharing the identity element $ e $. Therefore, the structure of the power graph of the dihedral group $ D_{2n} $ can be obtained from the power graph $ \mathcal{P}(\mathbb{Z}_{n}) $ by adding the $ n $ pendent vertices at the identity vertex $ e .$
If $ n=p^z $, where $ z $ is positive integer, then it is easy to see that
\[  \mathcal{P}(D_{2n}) = P_{3}[K_{n-1}, K_{1}, \overline{K}_{n}], \]
that is, the power graph $ \mathcal{P}(D_{2n}) $ of dihedral group for $ n=p^z $ is a \emph{pineapple} graph (graph obtained from the complete graph by attaching pendent vertices to any vertex of the complete graph). By using Theorem \ref{joined union}, $ \alpha n-1~\text{and} ~ \alpha $ are the $ A_{\alpha} $-eigenvalues of $ \mathcal{P}(D_{2n}) $ with multiplicities $ n-2 $ and $ n-1 $, respectively. The remaining three $ A_{\alpha} $-eigenvalues are given by the matrix
\[ \begin{pmatrix}
n-2+\alpha & 1-\alpha & 0\\
(1-\alpha)(n-2) & \alpha(2n-1) & (1-\alpha)n\\
0 & 1-\alpha & \alpha\\
\end{pmatrix}.
 \]
\qed

The following result gives the $ A_{\alpha} $-spectrum of the generalized quaternion group.
\begin{theorem}\label{dicyclic}
If $ n $ is a power of $ 2 $, then the $ A_{\alpha} $-spectrum of $ \mathcal{P}(Q_{n}) $ is 
\begin{equation*}
\{4 \alpha n-1,~ (2 \alpha n-1)^{(2n-3)}, ~ (4\alpha-1)^{(n)}, (1+2\alpha)^{(n-2)},~ x_{1},~ x_{2},~x_{3}\},
\end{equation*}
where $ x_{1},~ x_{2},~x_{3} $ are the zeros of the characteristic polynomial of the matrix given by \eqref{Qmat of quartenians 2}.
\end{theorem}
\textbf{Proof.} Consider the quaternion group $Q_{n}=<a,b ~|~ a^{2n}=e, b^2=a^n, ab=ba^{-1}>$. By definition of power graph, the identity element is always adjacent to every other vertex of $ \mathcal{P}(Q_{n}) $. In particular, if $ n $ is a power of $ 2 $, then by one of  the observations in  \cite{panda}, we found that $ a^n $ is also adjacent to all other vertices of $ \mathcal{P}(Q_{n}) $. Thus, two vertices of the power graph $  \mathcal{P}(Q_{n}) $ are connected to every other vertex of $  \mathcal{P}(Q_{n}) $. Since  $ n $ copies of $ K_{2} $ share the identity, so $  \mathcal{P}(Q_{n}) $ can be written as \[ \mathcal{P}(Q_{n})=K_{1,n+1}[K_{2},K_{2n-2},\underbrace{K_{2},K_{2},\dots,K_{2}}_{n}]. \]
By using Theorem \ref{joined union}, we see that $ \alpha_1= 4n-2, \alpha_{2}= \alpha_3=\dots=\alpha_{n+2}= 2 $ and the $ A_{\alpha} $-spectrum of $ \mathcal{P}(Q_{n}) $ consists of the simple eigenvalue $ 4\alpha n -1$, the eigenvalue $ 2\alpha n-1 $ with multiplicity $ 2n-3 $, the eigenvalue $ 4\alpha-1 $ with multiplicity $ n $ and the remaining $ n+2 $ eigenvalues are given by the following equitable quotient matrix
\begin{equation}\label{Qmat of quartenians}
M(Q_{n})=\begin{pmatrix}
1+\alpha(4n-2) & (1-\alpha)(2n-2) & (1-\alpha)2 & \dots&(1-\alpha)2\\
(1-\alpha)2 & 2n-3+2\alpha & 0 & \dots& 0 \\
(1-\alpha)2 & 0 & 1+2\alpha & \dots &  0 \\
\vdots & \vdots & \vdots & \dots & \vdots \\
(1-\alpha)2 & 0 & 0 & \dots & 1+2\alpha
\end{pmatrix}.
\end{equation}
Now, applying Lemma \ref{Lemma p2} with $ X=\begin{pmatrix}
1+\alpha(4n-2) & (1-\alpha)(2n-2)\\
(1-\alpha)2 & 2n-3+2\alpha
\end{pmatrix}, ~ \beta=\begin{pmatrix}
2(1-\alpha)\\ 0
\end{pmatrix}, B=(1+2\alpha) $ and $ C=(0) $, we see that $ 1+2\alpha $  is the $ A_{\alpha} $ eigenvalue of the matrix \eqref{Qmat of quartenians}  with multiplicity $ n-1 $. The remaining three $ A_{\alpha} $ eigenvalues of the matrix $ M(Q_{n}) $ in Equation \eqref{Qmat of quartenians} are given by the following equitable  quotient matrix 
\begin{equation}\label{Qmat of quartenians 2}
\begin{pmatrix}
1+\alpha(4n-2) & (1-\alpha)(2n-2) & 2n(1-\alpha)\\
2(1-\alpha) & 2n-3+\alpha & 0\\
2(1-\alpha) & 0 & 1+2\alpha
\end{pmatrix}.
\end{equation}\qed

We raise the following questions about $A_{\alpha}$-spectrum of $ \mathcal{P}(D_{2n}) $ and $ \mathcal{P}(Q_{n}) $.

\begin{problem}
Find the $A_{\alpha}$-spectrum of $ \mathcal{P}(D_{2n}) $ for the values of $ n\in\{pq,pqr, p^2q, (pq)^2\} $ and generalize (if possible) for any $ n $.
\end{problem}
\begin{problem}Find the $A_{\alpha}$-spectrum of  $ \mathcal{P}(Q_{n}) $ for the values of $ n\in\{pq,pqr, p^2q, (pq)^2\} $ and generalize (if possible) for any $ n $.
\end{problem}

Since for $ \mathcal{G}=\mathbb{Z}_{n}, ~ n=p, pq $, where $ p<q $ are primes, the multiplicity of $ A_{\alpha} $ eigenvalue $ \alpha n-1 $ is exactly $ \phi(n) $. Thus the partial solution of Problem \ref{Problem 1} for power graph $ \mathcal{P}(\mathbb{Z}_{n}) $ are graphs $ \mathcal{P}(\mathbb{Z}_{p})  $ and $ \mathcal{P}(\mathbb{Z}_{pq}) $. In order to show that there no other power graphs for $ \mathbb{Z}_{n} $, we have to show that $ x-(\alpha n-1) $ is not the factor of the characteristic polynomial of matrix  \eqref{quotient matrix of Z_n}. Similarly, in case of generalized quaternion group $ Q_{n}, ~ \text{for}~ n=p^z, ~ z\in \mathbb{N} $, multiplicity of $ \alpha n-1 $ is $ 1=\phi(2) $ and equality does not hold for any other $ n $. But as a whole Problem \ref{Problem 1} is still open.

\end{document}